\documentclass[9pt]{amsart}
\usepackage{amsfonts}
\usepackage{graphicx}
\usepackage{amscd}
\usepackage{amsmath,amsfonts,amssymb,amsthm,mathrsfs,xypic}
\newtheorem{thm}{Theorem}[section]

\newtheorem{lem}[thm]{Lemma}

\newtheorem{prop}[thm]{Proposition}

\newtheorem{Def}[thm]{Definition}
\newtheorem{ex}[thm]{\bf Example}
\numberwithin{equation}{section}

\def\ot{\otimes}

\newcommand{\Hom}{\operatorname{Hom}}

\begin{document}
\title [Gorenstein homological invariant properties under Frobenius extensions]
{Gorenstein homological invariant properties under Frobenius extensions}
\author [Zhibing  Zhao]
{Zhibing Zhao}
\thanks{Supported by Natural Science Foundation of China (No.11571329), the Natural Science Foundation of Anhui Province (No.1708085MA01) and Project of University Natural Science Research of Anhui Province (No.KJ2015A101).}
\thanks {Email: zbzhao$\symbol{64}$ahu.edu.cn}\maketitle
\begin{center}
School of Mathematical Sciences, Anhui University, Hefei, Anhui 230039
\end{center}
\begin{abstract} \ We prove that for a Frobenius extension, a module over the extension ring is Gorenstein projective if and only if its underlying
 module over the base ring is Gorenstein projective. For a separable Frobenius extension between Artin algebras, we obtain that
  the extension algebra is CM-finite (resp. CM-free) if and only if so is the base algebra. Furthermore, we prove that the reprensentation dimension of Artin algebras is invariant under separable Frobenius extensions.
 \vskip10pt

 \noindent 2010 Mathematics Subject Classfication: 13B02, 16E10, 16G10, 16G50

\noindent Key words: Frobenius extensions, separable extensions, Gorenstein projective modules, reprensentation dimension
\end{abstract}
\section{Introduction}
 Frobenius extensions are firstly introduced by Kasch in \cite{Kas1} as a generalization of Frobenius algebra. The fundamental example is the group
 algebras induced by a finite index subgroup. After that, some natural generalizations of different kinds of Frobenius extensions are defined
by  Nakayama-Tzuzuku, M$\ddot{u}$ller and Morita \cite{M,MU,NT}. There are other examples of Frobenius extensions including Hopf subalgebras, finite extensions of
enveloping algebras of Lie super-algebra and finite extensions of enveloping algebras of Lie coloralgebras  etc \cite{FMS,Sch}.

 Separable extensions are defined firstly by Hirata and Sugano in \cite{Hi} as a generalization of separable algebras, and they made a thorough
 study of these connection with Galois theory for
 noncommutative rings and generalizations of Azumaya algebras.  Separable extensions are closely related to Frobenius extensions. In fact, Frobenius extensions
 is a generalization of Frobenius algebras, while a separable algebra is a Frobenius algebra if it is finite projective over the base ring. A ring
  extension that is both a separable extension and a Frobenius extension is called a {\it separable Frobenius extension}. Sugano proved that the central projective
  separable extensions are Frobenius
 extensions in \cite{S}. More examples of separable Frobenius extensions can be found in Example \ref{examles}.
  We refer to a lecture due to Kadison \cite{Kad} for more details.

Recall from \cite{EO1} that an $A$-module $G$ is called Gorenstein projective if there is an exact sequence of projective $A$-modules
$$\mathbf{P}:=\cdots\rightarrow P_1\rightarrow P_0\rightarrow P^0\rightarrow P^1\rightarrow\cdots$$
in Mod-$A$ which is still exact after applying the functor ${\rm Hom}_A(-,Q)$ for any projective $A$-module $Q$ and $G={\rm Ker}(P_0\rightarrow P^0)$. In this case, the sequence $\mathbf{P}$ is called a {\it complete projective resolution} of $G$.
It is well-known that the notion of Gorenstein projective modules is a natural extension of  totally reflexive modules
to unnecessarily finitely generated modules (See \cite{EO1}) .
The subcategory of all Gorenstein projective $A$-modules (resp. all finitely generated Gorenstein projective $A$-modules) is denoted by $\mathcal{GP}(A)$ (resp. $\mathcal{FGP}(A)$).
The Gorenstein projective dimension of an $A$-module $M$, denoted by ${\rm Gpd}(M)$, is defined as
 ${\rm Gpd}(M)={\rm Inf}\{n\mid \exists$ an exact sequence $ 0\rightarrow G_n\rightarrow \cdots G_1\rightarrow G_0\rightarrow M\rightarrow 0$ with $G_i$ Gorenstein projective$\}$. We set
 ${\rm Gpd}(M)=\infty$ if there is no such a Gorenstein projective resolution. As a hot topic of Gorenstein homological algebra, Gorenstein projectivity of modules and
 Gorenstein projective dimension of modules are studied widely (see \cite{CH1,Ho,LZ,R}).

 The first motivation of this paper is inspired by the question in \cite{CH}. In \cite{CH}, Chen introduced a generalization of Frobenius extensions, called the totally reflexive extension,
 and he proved that totally reflexivity of modules is preserved under this kind of ring extension. At the end of \cite{CH}, Chen pointed out a question:
 is this true for Gorenstein projectivity in general?

 In \cite{RW}, Ren gave a partial answer for the above question. He proved that for a Frobenius extension, if a module over an extension ring is Gorenstein projective then
 it is also Gorenstein projective as a module over the base ring. Futhermore, the converse holds if the ring extension is a separable extension at the same time (see \cite [Theorem A] {RW}).
 The first main result in this paper extends the result in \cite{RW}. See Theorem \ref{GP property under FE}.

\noindent{\bf Theorem A.} {\it Let  $A/S$ be a Frobenius extension of rings and $M$ a right $A$-module. Then $M_A$ is a Gorenstein projective module if and
only if $M_S$ is a Gorenstein projective module.}

Furthermore, we prove the following statements.

\noindent{\bf Theorem B.} {\it  Let $A$ be an algebra over a commutative Artin ring $S$. If $A/S$ is a separable Frobenius extension, then we have

$(1)$ $A$ is CM-finite if and only if so is $S$;

$(2)$ $A$ is CM-free if and only if so is $S$.}

It is well-known that determining the representation type of an algebra is fundamental and important in representation theory of Artin algebra. Auslander proved that
there exists a 1-1 correspondence between the Morita equivalent classes of Artin algebras of finite representation type and that of Artin algebras with global dimension
at most 2 and with dominant dimension at least 2 (see \cite{A}). Motivated by this correspondence, Auslander introduced the notion of the representation dimension of Artin algebras. He
proved that an Artin algebra is of finite representation type if and only if its representation dimension is at most 2. In this sense, the representation dimension of an
Artin algebra is regarded as a trial to give a reasonable way of measuring homologically how far an Artin algebra is from being of finite representation type. Guo
 proved in \cite{G} that the representation dimension of an Artin algebra is invariant under stable equivalences.

 The second motivation of this paper comes from the results in \cite{Hu}. In \cite{Hu}, Huang and Sun proved that the representation dimension of Artin algebras is preserved under excellent
 extensions. It follows from \cite [Lemma 4.7]{Hu} that an excellent extension is a Frobenius extension. And an excellent extension is a separable Frobenius extension if the extension
 ring is commutative. But a separable Frobenius extension is not an excellent extension in general,
 see Example \ref{examles} for more details.
 As a separable Frobenius extension version of Theorem 4.8 in \cite{Hu}, we prove the following

\noindent{\bf Theorem C.} {\it Let $S$ be a commutative Artin ring and $A$ an $S$-algebra. If $A$ is a separable Frobenius extension of $S$, then
${\rm rep.dim}(A)={\rm rep.dim}(S)$.}

The paper is organized as follows. In Section 2, we give some notations in our terminology and some preliminary results which are used in this paper. In Section 3, we
prove the main result: Theorem A and B, see Theorem \ref{GP property under FE} and \ref{CM-finite and CM-free under SE} respectively.
In Section 4, the representation dimension under separable Frobenius extension of Artin algebras are studied, the main result Theorem C is proved.

Throughout this paper, all rings are associative rings with identity and all modules are unital right modules unless stated otherwise. Let $A$ be a ring.
We denote the category of all right $A$-module (resp. finitely generated right $A$-module) by Mod-$A$ (resp. mod-$A$).
\section{Preliminaries}

A {\it ring extension} $A/S$ is a ring homomorphism $\xymatrix@C=0.5cm{
   S \ar[r]^{l} &  A. }$ A ring extension is an algebra if $S$ is commutative and $l$ factor as the composition $S\rightarrow Z(A)\hookrightarrow A$ where
   $Z(A)$ is the center of $A$. The natural bimodule structure of ${_SA_S}$ is given by $s\cdot a\cdot s^\prime:=l(s)\cdot a\cdot l(s^\prime)$.
   Similarly, we can get some other module structures, for example $A_S$, ${_SA_A}$ and ${_AA_S}$, etc.

For any ring extension $A/S$, there is a {\it restriction functor} $R: {\rm Mod}$-$A\rightarrow {\rm Mod}$-$S$ which sends $M_A\mapsto M_S$, given by $m\cdot s:=m\cdot l(s)$.
In the opposite direction, there are two natural functors as follows:

$(1)$ $T=-\ot_SA_A :{\rm Mod}$-$S\rightarrow {\rm Mod}$-$A$ is given by $M_S\mapsto M\ot_S A_A$.

$(2)$ $H={\rm Hom}_S({_AA_S},-):{\rm Mod}$-$S\rightarrow {\rm Mod}$-$A$ is given by $M_S\mapsto {\rm Hom}_S({_AA_S},M_S)$.

It is easy to check that both $(T,R)$ and $(R,H)$ are adjoint pairs.
\begin{Def}$($\cite [Theorem 1.2]{Kad}$)$\label{Def of FE} A ring extension $A/S$ is a {\it Frobenius extension}, provided that one of the following equivalent conditions holds:

$(1)$ The functors $T$ and $H$ are naturally equivalent.

$(2)$ ${_SA_A}\cong {\rm Hom}_S({_AA_S},{_SS_S})$ and $A_S$ is finitely generated projective.

$(3)$ ${_AA_S}\cong {\rm Hom}_{S^{op}}({_SA_A},{_SS_S})$ and $_SA$ is finitely generated projective.

$(4)$ There exists $E\in {\rm Hom}_{S-S}(A,S)$, $x_i, y_i\in A$ such that $\forall a\in A$, one has $\sum\limits_ix_iE(y_ia)=a$
and $\sum\limits_iE(ax_i)y_i=a$.
\end{Def}

\begin{Def} A ring extension $A/S$ is a {\it separable extension} if and only if
$$\mu: A\ot_SA\rightarrow A, \quad a\otimes b\mapsto ab,$$
is a split epimorphism of $A$-$A$-bimodules. If a ring extension $A/S$ is both a Frobenius extension and a separable extension, then
it is called a {\it separable Frobenius extension}.
\end{Def}

Let $A/S$ be a ring extension and $M\in$Mod-$A$. Then $M_S$ is a right $S$-module. There is a natural surjective map
 $\pi: M\ot_SA\rightarrow M$ given by $m\ot a\mapsto ma $ for any $m\in M$ and $a\in A$. It is easy to check that $\pi$ is split as a homomorphism
of $S$-modules, we denoted by $M_S\mid M\ot_SA_S$ that $M$ is a direct summand of $M\ot_SA$ as right $S$-modules. However, $\pi$ is not split as an $A$-homomorphism in general.
 The following lemma comes from \cite{RW}, which is analogous to the results in \cite{p} for separable algebras over commutative rings.
\begin{lem}$\rm ($\cite [Lemma 2.9]{RW}$\rm )$\label{Property of SE} Let $A/S$ be a ring extension. The following are equivalent:

$(1)$ $A/S$ is a separable extension.

$(2)$ For any $A$-$A$-bimodule $M$, $M\ot_SA\rightarrow M$ is a split epimorphism of $A$-$A$-bimodules.

$(3)$ There exists an element $e\in A\ot_SA$, such that $\mu(e)=1_A$ and $ae=ea$ for any $a\in A$.
\end{lem}

\begin{ex}\label{examles}$(1)$ $($\cite [Example 2.10] {RW}$)$  For any finite group $G$, $\mathbb{Z}G/\mathbb{Z}$ is a separable Frobenius extension. Indeed,
 let $e=\frac{1}{\mid G\mid}\sum_{g\in G}g\ot_{\mathbb{Z}}g^{-1}\in \mathbb{Z}G\ot_{\mathbb{Z}}{\mathbb{Z}G}$, where $\mid G\mid$ is the order of $G$.
 Then $e$ satiefies the condition $(3)$ of the above lemma.

 $(2)$ $($\cite [Example 2.7 and 2.14]{Kad}$)$ Let $K$ be a field and $A=M_4(K)$. Let $S$ be the subalgebra of $A$ with $K$-basis consisting of the idempotents and matrix
 units,
 $$e_1=e_{11}+e_{14}, e_2=e_{22}+e_{33}, e_{21}, e_{31}, e_{41}, e_{42}, e_{43}.$$
 Then $A/S$ is a separable Frobenius extension. But since $A_S$ is not free as a right $S$-module, $A/S$ is not an excellent extension.

 $(3)$ Let $S$ be a commutative algebra and $A$ an {\it Azumaya algebra} over $S$. Then $A/S$ is a separable Forbenius extension.
  See \cite[Chapter 5]{Kad} for more details. We note that $A$ is not free over $S$, so $A/S$ is not an excellent extension in general.

 $(4)$ Every {\rm strongly separable extension} $($see \cite [Definition 3.1] {kad0}$)$ is a separable Frobenius extension. Specific examples of strongly separable extensions can be found in \cite{kad0}.
\end{ex}

 Let $\mathcal{C}$ be a subcategory of mod-$A$ and $M\in$mod-$A$. A homomorphism
 $f:C\rightarrow M$ in mod-$A$ is called a {\it right $\mathcal{C}$-approximation} of $M$ if $C\in\mathcal{C}$ and the sequence
 $ { \rm Hom}_A(-,C)\stackrel{(-,f)}\longrightarrow
{\rm Hom}_A(-,M)\longrightarrow 0$
is exact in $\mathcal{C}$. We say that an exact sequence
$0\longrightarrow
C_n\stackrel{f_n}\longrightarrow C_{n-1}
\stackrel{f_{n-1}}\longrightarrow \cdots\stackrel{f_1}\longrightarrow C_0\stackrel{f_0}\longrightarrow M\longrightarrow 0$
in mod-$A$ is an {\it $n$-$\mathcal{C}$-resolution} of $M$ if $C_i\in\mathcal{C}$ for any $0\leq i\leq n$, and the sequence
$$\xymatrix@C=0.5cm{
  0 \ar[r] & {\rm Hom}_A(-,C_n) \ar[rr]^{(-,f_n)} && {\rm Hom}_A(-,C_{n-1}) \ar[rr]^{(-,f_{n-1})} && \cdots}$$
  $$ \xymatrix@C=0.5cm{ \ar[rr]^{(-,f_1)} &&  {\rm Hom}_A(-,C_{0}) \ar[rr]^{(-,f_0)}
  && {\rm Hom}_A(-,M) \ar[r] & 0 }$$
is exact in $\mathcal{C}$ (See \cite{AR}).

We denoted by ${\rm add}M$ the full subcategory of mod-$A$ consisting of all modules isomorphic to direct summands of finite direct sums of copies of $M$,
and denoted by ${\rm Gen}M$ the full subcategory of mod-$A$ consisting of all modules $X$ such that there exists an epimorphism $M_0\twoheadrightarrow X$ with
$M_0\in {\rm add}M$. The following lemma comes from \cite{Hu}.

\begin{lem}$\rm ($\cite [Lemma 4.2]{Hu}$\rm )$\label{add-approximation} Let $A$ be an Artin algebra and $M, X\in {\rm mod}$-$A$. If $X=X_1\oplus X_2\in {\rm Gen}M$ has an $n$-${\rm add}M$-resolution:
$$\xymatrix@C=0.5cm{
     0 \ar[r] & M_n \ar[rr]^{f_n} && M_{n-1} \ar[rr]^{f_{n-1}} && \cdots \ar[rr]^{f_1} && M_0 \ar[rr]^{f_0} && X \ar[r] & 0 },$$
then $X_1$ has an $n$-${\rm add}M$-resolution:
$$\xymatrix@C=0.5cm{
     0 \ar[r] & M_n^\prime \ar[rr]^{f_n^\prime} && M_{n-1}^\prime \ar[rr]^{f_{n-1}^\prime} && \cdots \ar[rr]^{f_1^\prime} && M_0^\prime \ar[rr]^{f_0^\prime} && X_1 \ar[r] & 0 }.$$
\end{lem}

Let $A$ be an Artin algebra and $M\in$mod-$A$.
Recall that $M$ is called a {\it generator-cogenerator} for mod-$A$ if every indecomposable projective module and also every indecomposable injective module in mod-$A$
are in add$M$. The following lemma comes from \cite{EHIS}.

\begin{lem}\label{Eqcondition of gldim of end M}$\rm($\cite [Lemma 2.1]{EHIS}$\rm)$ Let $A$ be an Artin algebra and $M$ a generator-cogenerator for ${\rm mod}$-$A$.
Then the following statements are equivalent for any $n\geq 3$.

$(1)$ Any indecomposable module in ${\rm mod}$-$A$ has an $(n-2)$-${\rm add}M$-resolution;

$(2)$ ${\rm gl.dim End}(M_A)\leq n$.
\end{lem}

\section{Gorenstein projective dimensions under Frobenius extensions}
In this section, we will prove that for a Frobenius extension, a module over the extension ring is Gorenstein projective if and only if
its underlying module over the base ring is Gorenstien projective. Moreover, we obtain some homological properties,
including the Gorenstein global dimension of rings and the  CM-finiteness and CM-freeness of Artin algebras,
 are invariant under separable Frobenius extensions.

For a ring (or an algebra) $A$, we denoted by  $\mathcal{P}(A)$ the full subcategory of Mod-$A$ consisting of all
projective right $A$-modules. For a right $A$-module $M$, we denote
the projective dimension of $M$ by  ${\rm pd}(M_A)$. In order to prove that the Gorenstein projectivity of modules is
 preserved under Frobenius extensions, we need the following lemma.

\begin{lem}$\rm($\cite [Lemma 2.3]{RW}$\rm)$\label{GP-projective modules} Let $A/S$ be a Frobenius extension of rings and $M$ a right $A$-module.
If the underlying module $M_S$ is Gorenstein projective, then the following statements hold.

$(1)$ For any projective right $A$-module $P$ and any $i\geq 1$, ${\rm Ext}^i_A(M,P)=0$;

$(2)$ $M\ot_SA_A\cong {\rm Hom}_S({_AA},M)$ is a Gorenstein projective right $A$-module.
\end{lem}
\begin{thm}\label{GP property under FE}Let  $A/S$ be a Frobenius extension of rings and $M$ a right $A$-module. Then $M_A$ is a Gorenstein projective module if and
only if $M_S$ is a Gorenstein projective module.
\end{thm}

\noindent{\bf Proof.} The necessity holds by Lemma 2.2 in \cite{RW}. For the sake of completeness, we give the proof as follows.
  Assume that $M_A$ is Gorenstein projective as a right $A$-module. There exists a complete projective resolution
$\mathbf{P}:=\cdots\rightarrow P_1\rightarrow P_0\rightarrow P^0\rightarrow P^1\rightarrow\cdots$ in Mod-$A$ such that ${\rm Hom}_A(\mathbf{P},Q)$ is also exact
for any projective $A$-module $Q$ and $M={\rm Im}(P_0\rightarrow P^0)$. By the assumption, ${A_S}$ is finitely generated projective as a right $S$-module. Hence
$(P_i)_S\cong P_i\ot_AA_S$ is
projective as an $S$-module for any $i$. For any projective $S$-module $L$, we have ${\rm Hom}_S({_AA},L)\cong L\ot_SA_A$ is also projective as a right $A$-module.
Hence ${\rm Hom}_S(\mathbf{P},L)\cong{\rm Hom}_S(\mathbf{P}\ot_AA_S,L)\cong{\rm Hom}_A(\mathbf{P}, {\rm Hom}_S({_AA},L))$ is exact. So
$\mathbf{P}\cong \mathbf{P}\ot_AA_S$ in Mod-$S$ is a complete projective resolution of $M_S$ and therefore $M_S$ is a Gorenstein projective module.

Conversely, if $M_S$ is Gorenstein projective as a right $S$-module. Then we have ${\rm Ext}^i_A(M,P)$\\$=0$ for any projective $A$-module $P$ and any $i\geq 1$ and $M\ot_SA_A\cong {\rm Hom}_S({_AA},M)$ are Gorenstein projective as right $A$-modules by Lemma \ref{GP-projective modules}. We need to construct a complete projective
resolution of $M_A$ in Mod-$A$.

Since $M\ot_SA_A$ is Gorenstein projective, there exists a short exact sequence
$0\longrightarrow
 M\ot_SA_A \stackrel{f}\longrightarrow  P_0
\longrightarrow G\longrightarrow 0$
 in Mod-$A$ such  that $P_0\in \mathcal{P}(A)$ and $G\in \mathcal{GP}(A)$.
  Applying the restriction functor, we get the exact sequence $\xymatrix@C=0.5cm{
  0 \ar[r] & M\ot_SA_S}$ $\xymatrix@C=0.5cm{\ar[rr]^{f} && P_0 \ar[r] & G \ar[r] & 0 }$ in Mod-$S$ such  that $P_0\in \mathcal{P}(S)$
  and $G,  M\ot_SA\in \mathcal{GP}(S)$ by the necessity. Let $Q$ be any projective right $S$-module, and $g: M\ot_SA_S\rightarrow Q$ any $S$-homomorphism.
   $$\xymatrix{
  & Q  &  \\
  0\ar[r]& M\ot_SA \ar[u]_{g}\ar[r]^{f} &  P_0\ar@{..>}[ul]_{\exists h} \ar[r] & G \ar[r] & 0}$$
There exists an $S$-homomorphism $h:P_0\rightarrow Q$ such that $hf=g$ since $G$ is Gorenstein projective as an $S$-module.

  Note that there are two maps $i:M\rightarrow {\rm Hom}_S(A,M)$ given by $i(m)(a)=ma$ and $\pi:M\ot_SA\rightarrow M $ given by $\pi(m\ot a)=ma$ for any $m\in M$ and $a\in A$.
Since $M\ot_SA_A \cong {\rm Hom}_S({_AA},M)$, there is an $A$-monomorphism, which still denoted it by $i$, $i:M\hookrightarrow M\ot_SA$ and it is split
as a homomorphism of $S$-modules.
It follows that we get an exact sequence $\xymatrix@C=0.5cm{
  0 \ar[r] & M_A \ar[rr]^{fi} && P_0 \ar[r] & G_0 \ar[r] & 0 }$ in Mod-$A$ where $G_0={\rm Coker}(fi)$. Applying the restriction functor, we have an exact sequence
 $\xymatrix@C=0.5cm{
  0 \ar[r] & M_S \ar[rr]^{fi} && P_0 \ar[r] & G_0 \ar[r] & 0 }$ in Mod-$S$ such that $M_S\in \mathcal{GP}(S)$ (since $M_S\mid M\ot_SA_S$ and $M\ot_SA_S$ is Gorenstein projective) and $P_0\in\mathcal{P}(S)$.
 We claim that $G_0\in \mathcal{GP}(S)$.

 Let $Q$ be any projective $S$-module and $\alpha:M\rightarrow Q$ any $S$-homomorphism. Then $\alpha\pi$ is an $S$-homomorphism from $M\ot_SA$ to $Q$. Hence,
 there exists an $S$-homomorphism $\beta:P_0\rightarrow Q$ such that $\beta f=\alpha\pi$. And so $\beta (fi)=\alpha (\pi i)=\alpha$.
 $$\xymatrix{
  & Q  &  \\
  0\ar[r]& M \ar[u]_{\alpha}\ar[r]^{fi} &  P_0\ar@{..>}[ul]_{\exists \beta} \ar[r] & G_0 \ar[r] & 0 \\
  &{M\ot_SA} \ar[u]_{\pi} \ar[ur]_{f} &  }$$
 It follows from the exact sequence $$\xymatrix@C=0.5cm{
   0 \ar[r] & {\rm Hom}_S(G_0,Q) \ar[r] & {\rm Hom}_S(P_0,Q) \ar[r] & {\rm Hom}_S(M,Q)
   \ar[r] & {\rm Ext}_S^1(G_0,Q) \ar[r] & 0 }$$ that
   ${\rm Ext}_S^1(G_0,Q)=0$. Hence $G_0$ is Gorenstein projective as an $S$-module by \cite [Corollary 2.11]{Ho}.

 Note that $G_0$ is  Gorenstein projective as a right $S$-module, then ${\rm Ext}^i_A(G_0,P)=0$ for any projective $A$-module $P$ and any $i\geq 1$ and $G_0\ot_SA_A\cong {\rm Hom}_S({_AA},G_0)$ is  Gorenstein projective as a right $A$-module by Lemma \ref{GP-projective modules}. Repeating this process, we get an exact sequence
 $\mathbf{Q}:0\rightarrow M\rightarrow P_0\rightarrow P_1\rightarrow P_2\rightarrow\cdots $ in Mod-$A$ with $P_i$ projective and  ${\rm Hom}_A(\mathbf{Q},P)$ is also exact for any
 projective right $A$-module $P$. Linking up the projective resolution of $M_A$ and $\mathbf{Q}$, we get a complete projective resolution of $M_A$.
 The proof is completed.\hfill$\square$

A functor $F\colon\mathcal{C}\rightarrow\mathcal{D}$ is called a {\it Frobenius functor} if there exists a functor
$G\colon\mathcal{D}\rightarrow\mathcal{C}$ such that both $(F,G)$
and $(G,F)$ are adjoint pairs (see \cite{M}). By Definition \ref{Def of FE}, we know that the functors $-\ot_SA_A\cong {\rm Hom}_S({_AA_S,-})$ induced by
 the {\it Frobenius bimodule} (see \cite [Defintion 2.1] {Kad}) ${_AA_S}$ are  Frobenius functors.

\begin{prop}Let  $A/S$ be a Frobenius extension of rings. Then the Frobenius bimodule ${_AA_S}$ induces a Frobenius functor from
$\mathcal{GP}(A)$ to $\mathcal{GP}(S)$.
\end{prop}
\noindent{\bf Proof.} Put that $T=-\ot_SA_A:{\rm Mod}$-$S\rightarrow {\rm Mod}$-$A$ given by $M_S\mapsto M\ot_S A_A$. Then $T\cong H={\rm Hom}_S({_AA_S},-) $
is a Frobenius functor with the restriction functor $R$ as a left and right adjoint  at the same time.

By Lemma \ref{GP-projective modules} and Theorem \ref{GP property under FE}, we get $T\mid_{\mathcal{GP}(S)}\subseteq \mathcal{GP}(A)$
and $R\mid_{\mathcal{GP}(A)}\subseteq \mathcal{GP}(S)$, respectively. Hence $T$ is a Frobenius functor from
$\mathcal{GP}(A)$ to $\mathcal{GP}(S)$.\hfill$\square$

 \begin{prop}\label{Gp-dimension under FE}Let $A/S$ be a Frobenius extension of rings and $M$ a right $A$-module. Then ${\rm Gpd}(M_S)\leq {\rm Gpd}(M_A)$. Moreover, if  $A/S$ is separable, then
 ${\rm Gpd}(M_S)= {\rm Gpd}(M_A)$.
 \end{prop}
 \noindent{\bf Proof.} Without loss of the generality, we assume that ${\rm Gpd}(M_A)=n<\infty$. There is a Gorenstein projective resolution of $M_A$ as follows
 $$0\rightarrow G_n\rightarrow G_{n-1}\rightarrow\cdots \rightarrow G_1\rightarrow G_0\rightarrow M\rightarrow 0$$
 in Mod-$A$ with $G_i$ Gorenstein projective for any $0\leq i\leq n$. Applying the restriction functor, we get the following exact sequence
 $$0\rightarrow (G_n)_S\rightarrow (G_{n-1})_S\rightarrow\cdots \rightarrow (G_1)_S\rightarrow (G_0)_S\rightarrow M_S\rightarrow 0$$
in Mod-$S$ with $(G_i)_S\cong G_i\ot_AA_S$ Gorenstein projective as a right $S$-module for $0\leq i\leq n$
 by Theorem \ref{GP property under FE}. Hence ${\rm Gpd}(M_S)\leq n= {\rm Gpd}(M_A)$.

Furthermore, if $A/S$ is separable, then $M_A\mid M\ot_SA_A$ by Lemma \ref{Property of SE}. Without loss of the generality, we assume that ${\rm Gpd}(M_S)=m\leq\infty$. There exists
a Gorenstein projective resolution of $M_S$ as follows
$$0\rightarrow L_m\rightarrow L_{m-1}\rightarrow\cdots \rightarrow L_1\rightarrow L_0\rightarrow M\rightarrow 0$$
in Mod-$S$ such that $L_i$ is Gorenstein projective as a right $S$-module for any $0\leq i\leq m$. Applying the functor $-\ot_SA_A$, we get the following exact sequence
$$0\rightarrow L_m\ot_SA_A\rightarrow L_{m-1}\ot_SA_A\rightarrow\cdots \rightarrow L_1\ot_SA_A\rightarrow L_0\ot_SA_A\rightarrow M\ot_SA_A\rightarrow 0$$
in Mod-$A$ and $L_i\ot_SA_A$ is Gorenstein projective as a right $A$-module for $0\leq i\leq m$ by Lemma \ref{GP-projective modules}. And so
 ${\rm Gpd}(M\ot_SA_A)\leq m$. Hence, by the  \cite [Corollary 2.11]{Ho}, ${\rm Gpd}(M_A)\leq m= {\rm Gpd}(M_S)$.\hfill$\square$

 Recall that the right Gorenstein global dimension of a ring $A$, denoted by $r.{\rm Ggl.dim}(A)$, defined as $r.{\rm  Ggl.dim}(A)=$sup$\{ {\rm Gpd}(M_A)\mid$ M is a right $A$-module$\}$ (see \cite{BM}). It follows from Proposition \ref{Gp-dimension under FE} that the right Gorenstein global dimension of rings is preserved under  separable Frobenius extensions.

\begin{thm}Let $A/S$ be a separable Frobenius extension of rings. Then $r.{\rm Ggldim}(A)\\= r.{\rm Ggldim}(S)$.
\end{thm}
\noindent{\bf Proof.} It is easy to see $r.{\rm Ggl.dim}(A)\leq r.{\rm Ggl.dim}(S)$ by Proposition \ref{Gp-dimension under FE}.
Let $M_S$ be any $S$-module. Also by Proposition \ref{Gp-dimension under FE},
${\rm Gpd}(M\ot_SA_S)\leq {\rm Gpd}(M\ot_SA_A)$.  Since $M_S\mid M\ot_SA_S$, ${\rm Gpd}(M_S)\leq {\rm Gpd}(M\ot_SA_A)$ by \cite[Proposition 2.19]{Ho}.
It follows that $r.{\rm Ggl.dim}(S)\leq r.{\rm Ggl.dim}(A)$.\hfill$\square$

\vskip5pt

 Recall from \cite{B} that an Artin algebra $A$ is {\it Cohen-Macaulay finite}, or simply, {\it {\rm CM}-finite}, if $A$ has
only finitely many isomorphism classes of indecomposable finitely generated Gorenstein projective modules in mod-$A$. It is easy to see that $A$ is CM-finite if and only if
there exists a module $G\in $mod-$A$ such that $\mathcal{FGP}(A)={\rm add}G$. Clearly, $A$ is CM-finite if $A$ is of finite representation type. Some other examples of CM-finite algebra can
be found in \cite{R}.
 An Artin algebra $A$ is called {\it Cohen-Macaulay free}, or simply, {\it {\rm CM}-free}, if any Gorenstein projective module in mod-$A$ is projective. It is well-known that
 $\mathcal{GP}(A)=\mathcal{P}(A)$ if ${\rm gl.dim}(A)<\infty$, and $A$ is {\rm CM}-free if ${\rm gl.dim}(A)<\infty$.

 In order to  prove that the  {\rm CM}-finiteness and {\rm CM}-freeness of Artin algebras are preserved under separable
 Frobenius extensions, we need the following easy observation, which is maybe known.

 \begin{lem}\label{pdim of FE} Let $A$ be a separable Frobenius extension over a commutative Artin ring $S$ and $M$ a right $A$-module. Then we have

$(1)$ ${\rm pd}(M_A)={\rm pd}(M_S)$;

$(2)$ $r.{\rm gl.dim}(A)=r.{\rm gl.dim}(S)$.
\end{lem}
\noindent{\bf Proof.} $(1)$ Without loss of generality, we assume that ${\rm pd}(M_A)=n<\infty$ and $N_S$ is any right $S$-module.
Then ${\rm Ext}^{n+i}_S(M,N)\cong{\rm Ext}^{n+i}_S(M\ot_AA_S,N)\cong{\rm Ext}^{n+i}_A(M,{\rm Hom}_S({_AA_S},N))$\\$=0$ for any $i\geq 1$.
Hence ${\rm pd}(M_S)\leq n={\rm pd}(M_A)$.
 Conversely, we can assume that ${\rm pd}(M_S)=m<\infty$ and $W_A$ is any right $A$-module. Since $A/S$ is a separable extension ,
$M_A\mid (M\ot_SA)_A$ by Lemma \ref{Property of SE}. It follows from
${\rm Ext}^{m+i}_A(M\ot_SA_A,W)\cong {\rm Ext}^{m+i}_S(M_S,{\rm Hom}_A({_SA_A},W))=0 $ that ${\rm Ext}^{m+i}_A(M,W)=0$ for any $i\geq 1$. Thus
${\rm pd}(M_A)\leq m={\rm pd}(M_S)$.

$(2)$ It follows from $(1)$ that $r.{\rm gl.dim}(A)\leq r.{\rm gl.dim}(S)$. Conversely, let $X_S$ be any right $S$-module. Then $X_S\mid X\ot_SA_S$.
Since ${\rm pd}(X\ot_SA)_S={\rm pd}(X\ot_SA)_A$ by $(1)$, ${\rm pd}(X_S)\leq{\rm pd}(X\ot_SA)_A$. It follws that $r.{\rm gl.dim}(S)\leq r.{\rm gl.dim}(A)$.\hfill$\square$

\begin{thm}\label{CM-finite and CM-free under SE}Let $A$ be an algebra over a commutative Artin ring $S$. If $A/S$ is a separable Frobenius extension, then we have

$(1)$ $A$ is CM-finite if and only if so is $S$;

$(2)$ $A$ is CM-free if and only if so is $S$.
\end{thm}
\noindent{\bf Proof.} $(1)$ Assume that $S$ is {\rm CM}-finite, there exists a module $G\in$mod-$S$ such that $\mathcal{FGP}(S)={\rm add}G$.
We claim that $\mathcal{FGP}(A)={\rm add}(G\ot_SA)_A$. By Lemma \ref{GP-projective modules}, $(G\ot_SA)_A$ is Gorenstein projective as a right $A$-module,
and so ${\rm add}(G\ot_SA)_A\subseteq \mathcal{FGP}(A)$. Let $M_A$ be any indecomposable Gorenstein projective $A$-module in mod-$A$. By Proposition \ref{GP property under FE},
$M_S$ is also Gorenstein projective as a right $S$-module. There exists a positive integer $n$ such that $M_S\mid G_S^n$. Hence $M\ot_SA_A\mid (G\ot_SA_A)^n$. Since
$M_A\mid M\ot_SA_A$ by Lemma \ref{Property of SE}, $M_A\mid (G\ot_SA_A)^n$. Thus $\mathcal{FGP}(A)\subseteq {\rm add}(G\ot_SA)_A$.

Conversely, if $A$ is {\rm CM}-finite, then there exists a module $W\in$mod-$A$ such that $\mathcal{FGP}(A)={\rm add}W$. It suffices to prove that
$\mathcal{FGP}(S)={\rm add}W_S$. By Proposition \ref{GP property under FE}, $W_S$ is Gorenstien projective as a right $S$-module. Hence ${\rm add}W_S\subseteq \mathcal{FGP}(S)$.
Let $X_S$ be any indecomposable Gorenstein projective $S$-module in mod-$S$. By Lemma \ref{GP-projective modules}, $X\ot_SA_A$ is also Gorenstein projective as a right
$A$-module, and so $X\ot_SA_A\mid W_A^n$ for some positive integer $n$. Applying the restriction functor, we get $X\ot_SA_S\mid W_S^n$. It follows from $X_S\mid X\ot_SA_S$ that
$X_S\mid W_S^n$. Hence $\mathcal{FGP}(S)\subseteq {\rm add}T_S$.

$(2)$ Assume that $S$ is {\rm CM}-free and $G_A$ is any Gorenstein projective right $A$-module. By Proposition \ref{GP property under FE}, $G_S$ is
Gorenstein projective as a right $S$-module. Hence $G_S$ is projective by assumption. It follows from Lemma \ref{pdim of FE}$(1)$ that $G_A$ is also projective as
a right $A$-module. Therefore $A$ is {\rm CM}-free.

Conversely, if $A$ is {\rm CM}-free and $Q_S$ is any Gorenstein projective right $S$-module. By Lemma \ref{GP-projective modules}, $Q\ot_SA_A$ is Gorenstein projective as
a right $A$-module. Hence $Q\ot_SA_A$ is projective by assumption.  It follows from Lemma \ref{pdim of FE}$(1)$ that $Q\ot_SA_S$ is  projective as a right $S$-module.
Hence $Q_S$ is a projective right $S$-module as a direct summand of $Q\ot_SA_S$. So $S$ is {\rm CM}-free. \hfill$\square$

\section{Representation dimensions under separable Frobenius extensions}
In this section, we will show that the representation dimension of Artin algebras is invariant under separable Frobenius extensions.

 Recall that a module $M\in$mod-$A$ is called an additive generator for mod-$A$ if any indecomposable module in mod-$A$ is in add$M$. Obviously, an Artin algebra
$A$ is of finite representation type if and only if $A$ has an additive generator.

\begin{prop}\label{lemma of FRT}Let $A$ be an algebra over a commutative Artin ring $S$. If $A/S$ is a separable Frobenius extension, 
then $A$ is of finite representation type if and only if so is $S$.
\end{prop}
\noindent{\bf Proof.} We assume that $A$ is of finite representation type. Then there is an additive generator $M_A$ for mod-$A$. We claim that $M_S\cong M\ot_AA_S$
is an additive generator for mod-$S$. Let $X_S\in$mod-$S$ be any indecomposable right $S$-module. Then $X\ot_SA_A\mid M_A^n$ for some positive integer $n$. Applying the
restriction functor, we get $X\ot_SA_S\mid M_S^n$. It follows from $X_S\mid X\ot_SA_S$ that $X_S\mid M_S^n$. Thus, $M_S$ is an additive generator for mod-$S$.

Conversely, if $S$ is of finite representation type and $G_S$ is an additive generator for mod-$S$. We claim that $G\ot_SA_A$ is an additive generator for mod-$A$.
Let $N_A$ be any indecomposable right $A$-module in mod-$A$. Then $N_S\cong N\ot_AA_S\mid G_S^m$ for some positive integer $m$. And we get $N\ot_SA_A\mid (G\ot_SA)^m_A$.
It follows from $N_A\mid N\ot_SA_A$ that $N_A\mid (G\ot_SA)^m_A$. Therefore $G\ot_SA_A$ is an additive generator for mod-$A$. \hfill$\square$

As a homological dimension of measuring homologically how far an Artin algebra is from being of finite representation type, the representation dimension of an Artin algebra $A$,
denoted by ${\rm rep.dim}(A)$, was defined as follows in \cite{A}.

 \begin{Def} Let $A$ be an Artin algebra. The {\it representation dimension} ${\rm rep.dim}(A)$ of $A$ is defined as ${\rm inf}\{{\rm gl.dimEnd}(M_A)\mid M$ is
 a generator-cogenerator for {\rm mod}-$A\}$ if $A$ is non-semisimple; and ${\rm rep.dim}(A)=1$ if $A$ is semisimple.
 \end{Def}

 Let $A$, $S$ be Artin algebras and $A/S$ a separable Frobenius extention. By Lemma \ref{pdim of FE}(2), ${\rm rep.dim}(A)=1$ if and only if ${\rm rep.dim}(S)=1$.
And it follows from Proposition \ref{lemma of FRT} that ${\rm rep.dim}(A)\leq2$ if and only if ${\rm rep.dim}(S)\leq 2$. In general, we have the following main result of this section.
\begin{thm} Let $S$ be a commutative Artin ring and $A$ an $S$-algebra. If $A$ is a separable Frobenius extension of $S$, then
${\rm rep.dim}(A)={\rm rep.dim}(S)$.
\end{thm}

\noindent{\bf Proof.} By the definition of Frobenius extensions, $A$ is an Artin algebra. Then the assertion holds true provided either ${\rm rep.dim}(A)$
or ${\rm rep.dim}(S)$ is at most 2 by Lemma \ref{pdim of FE}(2) and Proposition \ref{lemma of FRT}.

Now assume that ${\rm rep.dim}(A)=n(\geq 3)$ and $M_A$ is a generator-cogenerator for mod-$A$ such that ${\rm gl.dimEnd}(M_A)=n$.
Let $X\in$mod-$S$ be indecomposable. It follows from $X_S\mid X\ot_SA_S$ that it is easy to see $M_S$ is a generator-cogenerator for mod-$S$.
Since $X\otimes_SA_A\in$mod-$A$, by Lemma \ref{Eqcondition of gldim of end M},
we have the following exact sequence
$$0\rightarrow M_{n-2}\rightarrow M_{n-3}\rightarrow \cdots\rightarrow M_0\rightarrow X\otimes_SA \rightarrow 0$$
in mod-$A$ (and hence in mod-$S$) with $M_i\in{\rm add}M_A$ such that
$$0\rightarrow {\rm Hom}_A(M,M_{n-2})\rightarrow {\rm Hom}_A(M,M_{n-2})\rightarrow\cdots $$
$$\rightarrow {\rm Hom}_A(M,M_0)\rightarrow {\rm Hom}_A(M, X\otimes_SA)\rightarrow 0$$
is also exact. By assumption, $A/S$ is a Frobenius extension, ${_SA}$ is finitely generated projective. Then
we get the following exact sequence
$$0\rightarrow {\rm Hom}_A(M,M_{n-2})\otimes_SA \rightarrow {\rm Hom}_A(M,M_{n-2})\otimes_SA \rightarrow\cdots$$
$$\rightarrow {\rm Hom}_A(M,M_0)\otimes_SA \rightarrow {\rm Hom}_A(M, X\otimes_SA)\otimes_SA \rightarrow 0.\eqno(*)$$
Since $({\rm Hom}_A({_SA},-)(\cong{-\ot_AA_S}), -\otimes_SA_A)$ is an adjoint pair, for any $N\in$mod-$A$ (also in mod-$S$) we get
 \begin{align*}{\rm Hom}_A(M,N)\ot_SA
&\cong {\rm Hom}_A(M, N\ot_SA)  \\
&\cong\Hom_S(\Hom_A(A,M),N) \\
&\cong\Hom_S(M,N),\end{align*}
where the first isomorphism comes from \cite [Theorem 3.2.14]{EO2} and the second isomorphism holds by the adjoint isomorphism.
Hence from the exact sequence $(*)$, we get the following exact sequence
$$0\rightarrow {\rm Hom}_S(M,M_{n-2})\rightarrow {\rm Hom}_S(M,M_{n-2})\rightarrow\cdots $$
$$\rightarrow {\rm Hom}_S(M,M_0)\rightarrow {\rm Hom}_S(M, X\otimes_SA)\rightarrow 0.$$
Thus $X\otimes_SA$ as an $S$-module has an $(n-2)$-add$M_S$-resolution. Since $X_S\mid (X\otimes_SA)_S$, $X_S$ has an
$(n-2)$-add$M_S$-resolution by Lemma \ref{add-approximation}. So ${\rm gl.dimEnd}(M_S)\leq n$ by Lemma \ref{Eqcondition of gldim of end M}
and therefore ${\rm rep.dim}(S)\leq n$.

Conversely, assume that ${\rm rep.dim}(S)=m(\geq 3)$ and $Q_S$ is a generator-cogenerator for mod-$S$ such that ${\rm gl.dimEnd}(Q_S)=m$.
Since $S\in{\rm add}Q_S$, $A_A\cong S\ot_SA_A\in {\rm add}(Q\ot_SA)_A$. It follows that $(Q\ot_SA)_A$ is a generator for mod-$A$.
On the other hand, if $Y\in$mod-$A$, then $Y$ is also a right $S$-module. Hence there exists a postive integer $t$ such that $0\rightarrow Y_S\rightarrow Q_S^{(t)}$ is exact
in mod-$S$, and so $0\rightarrow Y\ot_SA\rightarrow (Q\ot_SA)^{(t)}$ is exact in mod-$A$.
By the assumption, $A/S$ is a separable Frobenius extension. So $Y_A\mid (Y\ot_SA)_A$, and hence $(Q\ot_SA)_A$ is a cogenerator for mod-$A$. Thus
$(Q\ot_SA)_A$ is a generator-cogenerator for mod-$A$.

Let $V\in$mod-$A$ be indecomposable. By Lemma \ref{Eqcondition of gldim of end M}, $V_S$ has an $(m-2)$-add$Q_S$-resolution as an $S$-module
$$\xymatrix@C=0.5cm{
  0 \ar[r] & Q_{m-2} \ar[rr]^{f_{m-2}} && Q_{m-3} \ar[rr]^{f_{m-3}} && \cdots \ar[rr]^{f_1} && Q_0 \ar[rr]^{f_0} && V_S \ar[r] & 0 }.$$
We claim that the following sequence
$$0\rightarrow Q_{m-2}\ot_SA\rightarrow Q_{m-3}\ot_SA\rightarrow\cdots\rightarrow Q_0\ot_SA\rightarrow V\ot_SA\rightarrow 0.\eqno(**)$$
is an $(m-2)$-add$(Q\ot_SA)_A$-resolution of $(V\ot_SA)_A$.

It is easy to see that $Q_i\ot_SA_A\in$add$Q\ot_SA_A$. Let $K_i={\rm Ker}f_i$ for any $0\leq i\leq m-2$ and $K_{-1}=V$. Because ${_SA}$ is a finitely generated projective $S$-module,
 we have the following exact sequence
$$0\rightarrow K_i\ot_SA \rightarrow Q_i\ot_SA\rightarrow  K_{i-1}\ot_SA\rightarrow 0,$$
which is exact both as right $A$-modules and as right $S$-modules for any $0\leq i\leq m-2$. So the sequence $(**)$ is exact in mod-$A$. On the other hand, we have the following sequence
$$0\rightarrow {\rm Hom}_S(Q,K_i)\rightarrow {\rm Hom}_S(Q,Q_i)\rightarrow {\rm Hom}_S(Q,K_{i-1})\rightarrow 0$$
in mod-$S$, which induces the following exact sequence
$$0\rightarrow {\rm Hom}_S(Q,K_i)\ot_SA\rightarrow {\rm Hom}_S(Q,Q_i)\ot_SA\rightarrow {\rm Hom}_S(Q,K_{i-1})\ot_SA\rightarrow 0$$
for any $0\leq i\leq m-2$. By Lemma 3.2.4 in \cite{EO2}, ${\rm Hom}_S(Q,L)\ot_SA\cong {\rm Hom}_A(Q\ot_SA,L\ot_SA)$ for any $L\in$mod-$S$. Hence the sequence
$$0\rightarrow{\rm Hom}_A(Q\ot_SA,K_i\ot_SA)\rightarrow{\rm Hom}_A(Q\ot_SA,Q_i\ot_SA)$$
$$\rightarrow{\rm Hom}_A(Q\ot_SA,K_{i-1}\ot_SA)\rightarrow 0$$
is also exact for any $0\leq i\leq m-2$, which implies that the following sequence
$$0\rightarrow{\rm Hom}_A(Q\ot_SA,Q_{m-2}\ot_SA)\rightarrow{\rm Hom}_A(Q\ot_SA,Q_{m-3}\ot_SA)\rightarrow\cdots$$
$${\rm Hom}_A(Q\ot_SA,Q_0\ot_SA)\rightarrow{\rm Hom}_A(Q\ot_SA,V\ot_SA)\rightarrow 0$$
is exact. The claim is proved.

Notice that $V_A\mid (V\ot_SA)_A$, so $V_A$ has an $(m-2)$-add$(Q\ot_SA)_A$-resolution by Lemma \ref{add-approximation}. Thus
${\rm gl.dimEnd}((Q\ot_SA)_A)\leq m$ by Lemma \ref{Eqcondition of gldim of end M} and therefore ${\rm rep.dim}(A)\leq m$. The proof is finished.
\hfill$\square$

\noindent{\bf Acknowledgements}

The research was completed during the author's visit at the University of Washington. He would like to thank Professor James Zhang for his hospitality.
The author thanks the referees for the helpful comments and valuable suggestions.

\vspace{0.5cm}




\begin{thebibliography}{50}
\bibitem{A} M.Auslander, Representation dimension of Artin algebras, Queen Mary College Math. Notes, Queen Mary College, London, 1971.
\bibitem{AR} M.Auslander and I.Reiten, Applications of contravariantly finite subcategories, Adv. Math. 86(1991) 111-152.
\bibitem{B} A.Beligiannis, Cohen-Macaulay modules, (co)torsion pairs and virtually Gorenstein algebras, J. Algebra 288 (2005) 137-211.
\bibitem{CH} X.W. Chen, Totally reflexive extensions and modules, J. Algebra 379 (2013) 322-332.
\bibitem{CH1} X.W.Chen, An Auslander-type result for Gorenstein-projective modules, Adv. Math. 218(2008) 2043-2050.
\bibitem{BM} D.Bennis and N.Mahdou, Gorensten global dimension, Pro. Amer. Math. Soc. 138(2010) 461-465.
\bibitem{EO1} E.E.Enochs and O.M.G.Jenda, Gorenstein injective and projective modules, Math. Z.220(1995) 611-633.
\bibitem{EO2}  E.E.Enochs and O.M.G.Jenda, Relative homological algebra, de Gruyter Exp. Math., vol. 30, Walter de Gruyter, Berlin, New York,
2000.
\bibitem{EHIS} K.Ermann, T.Holm, O.Iyama and J.Schr\"{o}er, Radical embeddings and representation dimension, Adv. Math. 185(2004) 159-177.
\bibitem{FMS} D.Fischman, S.Montgomery and H.-J.Schneider, Frobenius extensions of subalgebras of Hopf algebras, Trans. Amer. Math. Soc. 349(1997) 4857-4895.
\bibitem{G} X.Guo, Representation dimension: An invariant under stable equivalence, Trans. Amer. Math. Soc. 357(2005)3255-3263.
\bibitem{Hi} K.Hirata and K. Sugano, On semisimple extensions and separable extensions over noncommutative rings, J. Math. Soc. Japan 18(1966) 360-373.
\bibitem{Ho} H.Holm, Gorenstein homological dimensions, J. Pure Appl. Algebra 189(2004), 167-193.
\bibitem{Hu} Z.Y.Huang and J. X. Sun, Invariant properties of represenations under excellent extensions, J. Algebra 358(2012) 87-101.
\bibitem{kad0} L.Kadison, The jones polynomial and certain separable Frobenius extensions, J. Algebra 186(1996) 461-475.
\bibitem{Kad} L.Kadison, New example of Frobenius extension, University Lecture Series, Vol 14, AMS. Provedence, Rhode Island, 1999.
\bibitem{Kas1} F.Kasch, Grundlagen einer theorie der Frobenius-Erweiterungen, Math. Ann., 127(1954) 453-474.
\bibitem{Kas2} F.Kasch, Projektive Frobenius-Erweiterungen,
Sizungsber. Heidelberger Akad. Wiss. Math.-Natur. Kl. (1960/61) 89-109.
\bibitem{LZ} Z.W.Li and P.Zhang, A construction of Gorenstein-projective modules, J. Algebra 323(2010) 1802-1812.
\bibitem{M} K.Morita, Adojint pairs of functors and Frobenius
extension, Sci. Rep. Toyko Kyoiku Daigaku (Sect. A) 9(1965) 40-71.
\bibitem{MU} B.M$\ddot{u}$ller, Quasi-Frobenius Erweiterungen I, Math. Zeit. 85(1964) 345-368.
\bibitem{NT} T.Nakamaya and T.Tsuzuku, On Frobenius extension I, Nagoya Math. J. 17(1960) 89-110.
\bibitem{p} R.S.Pierce, Associative algebra, GTM 88, New York: Springer, 1982.
\bibitem{R}  C.M.Ringel,  The Gorenstein-projective modules for the Nakayama algebras I, J. Algebra 385(2013) 241-261.
\bibitem{RW} W.Ren, Gorenstein projective mdoules and Frobenius extensions, Sci. China Ser. A. to apper, arXiv:1707.05885.
\bibitem{Sch} H.-J.Schneider, Normal basis and transitivity of crossed products for Hopf algebras, J. Algebra 151(1992) 289-312.
\bibitem{S} K.Sugano, Separable extensions and Frobenius extensions, Osaka J. Math. 7 (1970) 291-299.

\end{thebibliography}
\end{document}